\documentclass[12pt,leqno]{article}
\usepackage[pdftex]{graphicx}
\usepackage{amsfonts}
\usepackage{amsmath}
\usepackage{amsmath, amscd, amsfonts, amsthm, tikz, times}
\usepackage{setspace}
\usepackage{amssymb,textcomp}

\newcounter{conjecture}
\newcounter{remark}

\newtheorem{theorem}{Theorem}
\newtheorem{lemma}{Lemma}

\newtheorem{corollary}{Corollary}

\def \be{\begin{equation}}
\def \ee{\end{equation}}
\def \bt{\begin{theorem}}
\def \et{\end{theorem}}
\def \bc{\begin{corollary}}
\def \ec{\end{corollary}}
\def \bea{\begin{eqnarray}}
\def \eea{\end{eqnarray}}
\def \bas{\begin{eqnarray*}}
\def \eas{\end{eqnarray*}}
\def \bl{\begin{lemma}}
\def \el{\end{lemma}}

\def \Ga{\Gamma}

\def \la{\lambda}

\def \vski{\vspace{12pt}}
\def \ff{\infty}

\def \CC{{\cal C}}
\def \DD{\Delta}

\def \LL{{\cal L}}

\def \({\left(}
\def \){\right)}

\def \bc{\begin{center} }
\def \ec{\end{center} }
\def \bs{\begin{slide} }
\def \es{\end{slide} }

\def\square{{\vcenter{\vbox{\hrule height.3pt
         \hbox{\vrule width.3pt height5pt \kern5pt
            \vrule width.3pt}
         \hrule height.3pt}}}}
\def\qed{{\hfill $\square$ \bigskip}}

\newcounter{cccases}
\begin{document}

\title{On the nonexistence of pseudo-generalized quadrangles}

\author{Ivan Guo$^{\,\rm 2}$, Jack H. Koolen$^{\,\rm 1}$, ~Greg Markowsky$^{\,\rm 2}$, and Jongyook Park$^{\,\rm 3}${\footnote{Corresponding author.}}\\
{\small {\tt ivan.guo@monash.edu} ~~ {\tt koolen@ustc.edu.cn} ~~ {\tt gmarkowsky@gmail.com} ~~ {\tt jongyook@knu.ac.kr}}\\
{\footnotesize{$^{\rm 1}$Wen-Tsun Wu Key Laboratory of the CAS and School of Mathematical Sciences, USTC, Hefei, China}}\\
{\footnotesize{$^{\rm 2}$School of Mathematics,  Monash University, Australia}}\\
{\footnotesize{$^{\rm 3}$Department of Mathematics, Kyungpook National University, Daegu, 41566, Republic of Korea}}}





\maketitle


\begin{abstract}
In this paper we consider the question of when a strongly regular graph with parameters $((s+1)(st+1),s(t+1),s-1,t+1)$ can exist. These parameters arise when the graph is derived from a generalized quadrangle, but there are other examples which do not arise in this manner, and we term these {\it pseudo-generalized quadrangles}. If the graph is a generalized quadrangle then $t \leq s^2$ and $s \leq t^2$, while for pseudo-generalized quadrangles we still have the former bound but not the latter. Previously, Neumaier has proved a bound for $s$ which is cubic in $t$, but we improve this to one which is quadratic. The proof involves a careful analysis of cliques and cocliques in the graph. This improved bound eliminates many potential parameter sets which were otherwise feasible.
\end{abstract}

\section{Introduction and statement of results}

Generalized quadrangles are a class of intensely-studied incidence structures consisting of a finite set of elements, which we may consider points, together with a collection of subsets of these points, which we consider lines, satisfying the following conditions.

\begin{itemize} \label{}

\item[(i)] There is an integer $s \geq 1$ such that each line contains exactly $s+1$ points. Any two distinct lines share at most one point.

\item[(ii)] There is an integer $t \geq 1$ such that each point is contained on exactly $t+1$ lines. Any two distinct points are contained on at most one line.

\item[(iii)] For every point $p$ not on a line $L$, there is a unique line $M$ and a unique point $q$, such that $p$ is on $M$, and $q$ is on $M$ and $L$.

\end{itemize}

For background and applications of generalized quadrangles, see \cite{payne}. We denote a generalized quadrangle with parameters $s$ and $t$ as $GQ(s,t)$. The collinearity graph $G$ of a generalized quadrangle with parameters $s$ and $t$, which we will also denote $GQ(s,t)$ in a mild abuse of notation, is {\it strongly regular} with parameters $(v,k,\lambda,\mu)=((s+1)(st+1),s(t+1),s-1,t+1)$; this means that $G$ is a $k$-regular graph on $v$ vertices, with the property that every pair of adjacent vertices have exactly $\lambda$ common neighbors and every pair of non-adjacent vertices have exactly $\mu$ common neighbors. We will refer to such a strongly regular graph as $srg((s+1)(st+1),s(t+1),s-1,t+1)$; see \cite{godroy} or \cite{cameron} for more on strongly regular graphs. The cases $t=1$ and $s=1$, corresponding to the $(s+1) \times (s+1)$ rook's graph and the complete bipartite graph $K_{t+1,t+1}$, respectively, are considered trivial, so we will assume $s,t \geq 2$ in what follows. A strong necessary condition for a $GQ(s,t)$ to exist is that $t \leq s^2$. If a $GQ(s,t)$ exists for some $s,t$, let us call it ${\cal Q}$, then we may obtain the dual structure ${\cal Q'}$ which has the lines and points of ${\cal Q}$ as its points and lines, respectively, with the concepts of colinearity of points and intersections of lines also being interchanged, and ${\cal Q'}$ is a $GQ(t,s)$. Thus, by this process of dualization, we obtain that $s \leq t^2$, as well.

\vski

On the other hand, it is entirely possible for a $srg((s+1)(st+1),s(t+1),s-1,t+1)$ to exist which is not isomorphic to the colinearity graph of a generalized quadrangle; a good example of this is the Cameron graph, which is a $srg(231,30,9,3)$ (see \cite{broupage} or \cite{camgraph}) and yet which is clearly not a generalized quadrangle since $s > t^2$. We refer to such a strongly regular graph, which has parameters $((s+1)(st+1),s(t+1),s-1,t+1)$ but which is not isomorphic to a $GQ(s,t)$, as a {\it pseudo-generalized quadrangle}, and denote it as $PGQ(s,t)$. Thus, the Cameron graph is a $PGQ(10,2)$. It is an interesting fact that the relation $t \leq s^2$ persists for any $PGQ(s,t)$ as a consequence of the Krein eigenvalue bounds (see the proof of \cite[Lemma 10.8.3]{godroy}), but it is evident that $s \leq t^2$ does not, and we have no concept of duality to apply to pseudo-generalized quadrangles.

\vski

Previously, the strongest known upper bound on $s$ in terms of $t$ was given by the following result, which is Neumaier's claw bound applied to the parameter set of a pseudo-generalized quadrangle.

\begin{theorem} \cite[Thm. 4.7 (iii)]{neum} \label{oldbusted}

$$s \leq \frac{t(t+1)(t+2)}{2}$$

\end{theorem}

See also \cite[Sec. 8.6]{broupage} for more on this and related results. We will provide the following improvement.

\begin{theorem} \label{newhotness}
Let $\theta$ and $\beta$ be integers satisfying $\theta \geq t+2$ and $2\leq \beta\leq t+1$, then
\[
s \leq \max\left\{\frac{t}{\theta-t}\binom{\theta+1}{2},\ t(2\theta-1),\ \binom{\beta}{2}t,\ (t+1)^2\theta\binom{\beta}{2}^{-1}\right\}.
\]

\end{theorem}

By choosing the parameters $\theta$ and $\beta$ appropriately, we obtain the following corollary.

\begin{corollary} \label{iu}
We have
\[
s \leq t\left\lfloor \frac{8}{3}t+1\right\rfloor
\]
for all pseudo-generalized quadrangles $PGQ(s,t)$ with $t\geq 2$.
\end{corollary}



%
%
%
%
%
%
%
%

%

In the next section we prove these results, while in the final section we list a number of otherwise feasible parameter sets which are eliminated by our result.

\section{Proofs}

{\bf Proof of Theorem \ref{newhotness}}

Let $G$ be a $PGQ(s,t)$. We will call the complete bipartite graph $K_{1,r}$ the {\it $r$-claw}, the vertex of degree $r$ will be the {\it center} of the $r$-claw, and the other vertices will be the {\it leaves}. We will use the standard notation $V(G)$ for the vertex set of $G$. For any $x \in V(G)$ let $\Ga(x) = \{y \in V(G): y \sim x\}$, and for any vertices $x_1, \ldots, x_m$ of $G$ with $m \geq 2$, let $\Ga(x_1, \ldots, x_m) = \cap_{j=1}^m \Ga(x_j)$. Let $\DD(x)$ denote the local graph at $x$, that is, $\DD(x)$ is the subgraph of $G$ induced on $\Ga(x)$. For any vertex $x$ in $V(G)$, define its \emph{claw number} $\phi(x)$ by

$$
\phi(x) = \max\{r: x \mbox{ is the center of an induced $r$-claw in } G\}.
$$

Equivalently, $\phi(x)$ is the order of the maximal induced coclique in $\DD(x)$. Let $d(x,y)$ denote the standard shortest-path distance between $x$ and $y$. We will prove Theorem \ref{newhotness} in a series of lemmas.

\begin{lemma} \label{lowerboundonphi}
For all $x \in V(G)$, $\phi(x) \geq t+1$.
\end{lemma}

{\bf Proof:} Suppose there is an $x$ with $\phi(x) = r \leq t$. Let $G'$ denote an induced $r$-claw with center at $x$, and $L$ be the set of leaves of $G'$. Every vertex in $\Ga(x)\setminus L$ must be adjacent to a point in $L$, since otherwise $G'$ could be extended to a larger induced claw, and thus

$$
k=|\Ga(x)| \leq |L|+ \sum_{y \in L} |\Ga(x,y)| \leq r+ r \la \leq t + t(s-1) =st.
$$

However, $k=s(t+1)$, and this contradiction proves the lemma. \qed

\begin{lemma} \label{oneguy}
If $\phi(x) = t+1$ for some $x \in V(G)$, then $\DD(x)$ is the disjoint union of $t+1$ cliques, each of which is maximal (i.e. can't be extended to a larger clique) and of order $s$.
\end{lemma}

{\bf Proof:} Choose $y_1, \ldots, y_{t+1}$ which induce a coclique in $\DD(x)$. The fact that $\phi(x) = t+1$ implies that the subgraph induced on $x,y_1, \ldots, y_{t+1}$ cannot be extended to a larger claw, and we conclude that $\Ga(x) = \cup_{j=1}^{t+1}(\{y_j\} \cup \Ga(x,y_j))$. Furthermore, since $|\Ga(x)| = k = s(t+1) = \sum_{j=1}^{t+1} |\{y_j\} \cup \Ga(x,y_j)|$ we conclude that $\Ga(x,y_i) \cap \Ga(x,y_j) = \emptyset$ for $i \neq j$. This implies that if we choose $y'_j \in \Ga(x,y_j)$ for some $j$, then $y_1, \ldots,y_{j-1},y'_j,y_{j+1}, \ldots, y_{t+1}$ also induce a coclique in $\DD(x)$, to which the same argument applies, and we see that $\{y_j\} \cup \Ga(x,y_j) = \{y'_j\} \cup\Ga(x,y'_j)$. It follows that each of the sets of the form $\{y_j\} \cup \Ga(x,y_j)$ induces a clique of order $s$ in $\DD(x)$, and therefore $\DD(x)$ is the disjoint union of $t+1$ cliques, each of order $s$. It is clear from the construction that these cliques are maximal. \qed

\begin{lemma} \label{itsgeo}
If $\phi(x) = t+1$ for all $x \in V(G)$, then in fact $G$ is a $GQ(s,t)$, and therefore $s \leq t^2$.
\end{lemma}

{\bf Remark:} This lemma can be deduced from the more general results in \cite{metsch}, but for the benefit of the reader we include the following simple proof, which is not as general as Metsch's arguments but which suffices for our purposes.

\vski

{\bf Proof:} Suppose $\phi(x) = t+1$ for all $x \in V(G)$. Lemma \ref{oneguy} shows that for every $x$ the local graph $\DD(x)$ is the disjoint union of cliques, and we will use these cliques to show that $G$ is a generalized quadrangle. Rephrasing Lemma \ref{oneguy}, each $x \in V(G)$ is a member of precisely $t+1$ cliques of order $s+1$, any two of which intersect only at $x$. Let $\LL$ be the set of all such cliques of order $s+1$ in $G$. We claim now that if $G'$ is an induced subgraph of $G$ which is a clique of order $r \geq 2$ then $G'$ is contained in a unique clique $C$ in $\LL$. In order to prove this, it clearly suffices to assume $r=2$, and $V(G') = \{x_1,x_2\}$. We may choose $C$ to be the clique formed by the clique-partition of $\DD(x_1)$ which contains $x_2$. If $x_1$ and $x_2$ are contained in any other clique $C'$ in $\LL$ then $C' \backslash \{x_1\}$ must also be part of the clique-partition of $\DD(x_1)$, and since $x_2 \in C$ and $x_2 \in C'$ we must have $C=C'$, proving the claim.

\vski

An incidence structure $GQ(s,t)$ can now be formed with the vertices of $G$ as points and $\LL$ as the set of lines. The claim we just proved shows that any two distinct lines in $\LL$ share at most one point, and that any two distinct points are contained on at most one line in $\LL$. Thus conditions $(i)$ and $(ii)$ from the introduction are satisfied, and we need only show $(iii)$.
Suppose $x$ is a vertex and $C$ is a clique in $\LL$ not containing $x$; we must show that there is precisely one clique in $\LL$ containing $x$ which has nonempty intersection with $C$, and that said intersection consists of a single point. We begin by noting that there is at least one point $z \in C$ such that $x \not \sim y$, since otherwise $C$ could be extended to a clique of order $s+1$. Let the cliques containing $z$ be labelled $D_1, \ldots , D_{t+1}$; $C$ must be equal to one of these. We claim now that $x$ must be adjacent to one point in each of the $D_j$'s; if not, since $\mu = t+1$, by the pigeonhole principle there must be vertices $w_1, w_2$ in some $D_j$ which are both adjacent to $x$. But then the subgraphs induced on $\{x,w_1,w_2\}$ and $\{z,w_1,w_2\}$ are both cliques, and both must therefore be contained in the unique clique containing $w_1$ and $w_2$. However, this implies that $z$ and $w$ are adjacent, a contradiction. This proves that $x$ must be adjacent to exactly one point in each of the $D_j$'s, in particular to one point $w \in C$. The edge $(x,w)$ can be extended to a unique clique $C'$, and $C$ and $C'$ intersect only at $w$. This proves $(iii)$, and the lemma. \qed

Let $\theta$ be an integer satisfying $\theta \geq t+2$.
\begin{lemma} \label{bound1}
If $\phi(x) \geq \theta+1$ for some $x \in V(G)$, then $s \leq \frac{t}{\theta-t}\binom{\theta+1}{2}$.
\end{lemma}

{\bf Proof:} The following result, translated into our notation, was proved in \cite{shilla}.

\begin{theorem} \cite{shilla} \label{shill}
If $G$ is a distance-regular graph, and if there is an $x \in V(G)$ which has $\phi(x) \geq r$, then

$$
\mu-1 \geq \frac{r(\la+1)-k}{{r \choose 2}}.
$$
\end{theorem}

Applying this in our situation ($\mu-1 = t, r = \theta+1, k=s(t+1), \la+1 = s$) and simplifying proves the lemma. \qed

This lemma proves Theorem \ref{newhotness} when $\phi(x) \geq \theta+1$ for some $x$. We are now reduced to the case $t+1\leq \phi(x) \leq \theta$ for all $x \in V(G)$, and furthermore $\phi(x) \geq t+2$ for at least one $x$.

\begin{lemma} \label{}
Fix $y \in V(G)$, and let $x_1, \ldots, x_{\phi(y)}$ be a set of vertices that induces a coclique in $\DD(y)$. Then each $x_j$ lies in an induced clique of order at least $s-t(\theta-1)$ in $\DD(y)$.
\end{lemma}

{\bf Proof:} For any $j \in \{1,\ldots, \phi(y)\}$, since $|\Ga(y,x_i,x_j)| \leq \mu - 1 = t$ for $i \neq j$, we have

\begin{equation} \label{} \nonumber
\begin{split}
|\{x_j\}\cup\Ga(y,x_j) \backslash \bigcup_{\stackrel{1 \leq i \leq \phi(y)}{i \neq j}}\Ga(y,x_i,x_j)| & \geq |\{x_j\}\cup\Ga(y,x_j)| - \sum_{\stackrel{1 \leq i \leq \phi(y)}{i \neq j}}|\Ga(y,x_i,x_j)| \\
& \geq s - t(\phi(y)-1).
\end{split}
\end{equation}

If $\{x_j\}\cup\Ga(y,x_j) \backslash \bigcup_{\stackrel{1 \leq i \leq \phi(y)}{i \neq j}}\Ga(y,x_i,x_j)$ is not a clique then it must contain two non-adjacent vertices, call them $z_1$ and $z_2$, and we see that the vertices $x_1, \ldots, x_{j-1}, z_1, z_2, x_{j+1}, \ldots , x_{\phi(y)}$ induce a coclique of order $\phi(y)+1$ in $\DD(y)$, a contradiction. The proof is completed by noting that $\phi(y) \leq \theta$. \qed

By extending each of the cliques in the previous lemma to a maximal clique, to any induced $\phi(y)$ claw with center $y$ and leaves $x_1, \ldots , x_{\phi(y)}$ we may associate maximal cliques $C_1, \ldots ,C_{\phi(y)}$ of order at least $(s+1) - t(\theta-1)$ where each $C_j$ contains $y$ and $x_j$.

\begin{lemma} \label{hungary}
If there exists a vertex $x_0$ of $\DD(y)$ which is not contained in any of the $C_j$'s, then $s \leq t(2\theta-1)$.
\end{lemma}

{\bf Proof:} Suppose that $s>t(2\theta-1)$. We will show that $\phi(y) \geq \theta+1$, which contradicts the assumption we are operating under. If $x_0 \in V(\DD(y))$ but $x_0 \notin C_j$ for all $j$, we begin by noting that, since each $C_j$ is maximal, each contains some vertex $z_j$ not adjacent to $x_0$, and it follows from this that $|C_j \cap \Ga(y,x_0)| \leq \mu-1 = t$. The same logic shows that if $z_i \in C_i$ for $i \neq j$ then $|C_j \cap \Ga(y,z_i)| \leq t$. This observation allows us to construct an induced $(\phi(y)+1)$-claw, as follows. Begin by choosing $x_1 \in C_1$ which is not in $\Ga(y,x_0)$; this is possible because $|C_1 \backslash \Ga(y,x_0)| \geq |(C_1 \backslash \{y\})| - |C_1 \cap \Ga(y,x_0)| \geq (s-t(\theta-1)) -t>0$. We then apply this argument inductively. Having chosen $x_j \in C_j$, we choose $x_{j+1} \in C_{j+1}$ to be a vertex in the set $C_{j+1} \backslash \cup_{i=0}^j \Ga(y,x_j)$, and again this is possible since $|C_{j+1} \backslash \cup_{i=0}^j \Ga(y,x_j)| \geq |C_{j+1}\backslash \{y\}| - \sum_{i=0}^{j}\Ga(y,x_j)| \geq (s-t(\theta-1)) -t(j+1) \geq (s-t(\theta-1)) -t\theta = s-t(2\theta-1) > 0$, by our assumption. By construction, the subgraph induced on $\{y, x_0, x_1, \ldots ,x_{\theta}\}$ is a $(\phi(y)+1)$-claw, and the result follows. \qed

So far, for every $y\in V(G)$, we have identified a sequence of maximal cliques $C_1,\ldots,C_{\phi(y)}$.
Let $\CC$ denote the set of all such maximal cliques, starting from all possible $y\in V(G)$.

\begin{lemma} \label{agent}
Suppose $s>t(2\theta-1)$. Then each $y \in V(G)$ is contained in precisely $\phi(y)$ cliques in $\CC$. Furthermore, if $x \sim y$ then there is a unique clique in $\CC$ which contains both $x$ and $y$.
\end{lemma}

{\bf Proof:} We will prove the second statement first. We note that Lemma \ref{hungary} shows that if $x \sim y$ then there is at least one clique in $\CC$ containing both of them. Suppose there are two, $C_1$ and $C_2$. If they do not coincide, then by maximality we can find $z_1 \in C_1, z_2 \in C_2$ with $d(z_1,z_2) = 2$. Then, since $(C_1 \cap C_2) \subseteq \Ga(z_1,z_2)$ we have $|C_1 \cap C_2| \leq \mu = t+1$. Thus, $|C_1 \cup C_2| = |C_1|+|C_2|-|C_1 \cap C_2| \geq 2(s+1-t(\theta-1)) - (t+1) = 2s+1 -t(2\theta-1)$. On the other hand, $(C_1 \cup C_2) \subseteq (\Ga(y,x) \cup \{x,y\})$ and therefore $|C_1 \cup C_2| \leq \la+2 = s+1$. We obtain $s+1 \geq 2s+1 -t(2\theta-1)$, hence $s \leq t(2\theta-1)$, a contradiction. As for the first statement, we know from earlier arguments that $y$ is contained in at least $\phi(y)$ cliques in $\CC$, each of which is associated to a leaf of an induced $\phi(y)$-claw with center at $y$. Furthermore, Lemma \ref{hungary} shows that the union of the vertices of these $\phi(y)$ cliques is $\{y\} \cup \Ga(y)$. From this and the second statement in this lemma we conclude that there are no cliques in $\CC$ other than these which contain $y$ (each vertex in $\DD(x)$ is associated to a unique clique in $\CC$), and therefore that $y$ is contained in precisely $\phi(y)$ cliques. \qed

Let $\tilde V \subseteq V(G)$ be the set of vertices whose claw numbers are at least $t+2$. In other words
\[
\tilde V=\{x\in V(G): \phi(x)\geq t+2\}.
\]

\begin{lemma} \label{steps}
Suppose $s >t(2\theta-1)$. Then $|\tilde V| \leq |\CC| \leq 2(st+1)\theta < |V(G)|$.
\end{lemma}

{\bf Proof:} Recall that $|V(G)|=(st+1)(s+1)$. If we multiply the number of vertices by the maximum number of cliques per vertex and then divide by the minimum number of vertices per clique we get an upper bound for $|\CC|$. This yields
\[
|\CC| \leq \frac{|V(G)|\theta}{1+s-t(\theta-1)} = \frac{(st+1)(s+1)\theta}{1+s-t(\theta-1)}.
\]
The assumption on $s$ implies easily that $s+1 > 2t(\theta-1)$. Hence $\frac{s+1}{1+s-t(\theta-1)}<2$ and $|\CC| \leq 2(st+1)\theta$.
Next, $ 2(st+1)\theta < |V(G)|$ follows from the fact that the assumption on $s$ implies $2\theta<s+1$.
\vski

Finally, it suffices to prove that $|\tilde V| \leq |\CC|$. Let us form the {\it vertex-clique incidence matrix} $R$ of $G$: this is an $n \times c$ matrix, where $n=|V(G)|$ and $c=|\CC|$, with entry $R_{ji} = 1$ if vertex $j \in V(G)$ is contained in clique $i \in \CC$ and $R_{ji} = 0$ otherwise. Lemma \ref{agent} shows that
$$
RR^T = A + D,
$$
where $A$ is the adjacency matrix of $G$ and $D$ is a diagonal matrix where $D_{jj}$ is equal to the claw number of vertex $j$. Since $|\CC|<|V(G)|$ and $RR^T$ is positive semidefinite, we must have $Rank(R) < n$ and the $n-Rank(R)$ smallest eigenvalues of $RR^T$ are zeros. From \cite[Lemma 10.82]{godroy}, the smallest eigenvalue of $A$ is $-(t+1)$. By the Courant-Weyl inequalities \cite[Theorem 2.8.1]{brouwer2011spectra}, this implies the $n-Rank(R)$ smallest eigenvalues of $D$ must be at most $t+1$.
Since $D$ is diagonal matrix with exactly $n-|\tilde V|$ diagonal elements being $t+1$ and $Rank(R)\leq c=|\CC|$, we can therefore conclude that $|\tilde V| \leq |\CC|$, completing the proof.
\qed

Let $\alpha$ be an integer satisfying $0\leq \alpha\leq t-1$. Let $\beta=t+1-\alpha$, so $2\leq \beta\leq t+1$.
We will consider two separate cases:

\vski

(A) There exists $y\in\tilde V$ such that $y$ belongs to at least $\alpha+1$ cliques of order $s+1$;\\
(B) For all $y\in\tilde V$, $y$ belongs to at most $\alpha$ cliques of order $s+1$.

\vski

Note that if a vertex belongs to a clique of order $s+1$, the clique must be maximal.

\begin{lemma} \label{case1}
Suppose $s >t(2\theta-1)$ and (A) holds, then $s\leq \binom{\beta}{2} t$.
\end{lemma}

{\bf Proof:} Suppose that $y\in\tilde V$ satisfies condition (A).
By the definition of $\tilde V$, $y$ is the center of a $\phi(y)$-claw where $\phi(y)= t+2+\gamma$ for some $\gamma\geq 0$. Call the leaves of this claw $x_1, \ldots , x_{\phi(y)}$. Without loss of generality, let $x_1, \ldots , x_{\alpha+1}$ belong to distinct cliques of order $s+1$. Let us count the vertices of $\DD(y)$ outside of these cliques.
By the inclusion-exclusion principle we have
\[
k-s(\alpha+1) = \bigg|\bigcup_{i=\alpha+2}^{\phi(y)} (\Ga(y,x_i) \cup \{x_i\})\bigg| \geq \sum_{i=\alpha+2}^{\phi(y)} |\Ga(y,x_i) \cup \{x_i\}| - \sum_{\alpha+2 \leq i < j \leq \phi(y)}|\Ga(y,x_i,x_j)|.
\]
The parameters of the graph give $k = s(t+1)$, $|\Ga(y,x_i) \cup \{x_i\}| = s$, and $|\Ga(y,x_i,x_j)| \leq t$ (because $x_i$ and $x_j$ are not adjacent, $\mu = t+1$, and $y \notin \Ga(y,x_i,x_j)$). Therefore, using $\phi(y)=t+2+\gamma=\alpha+1+\beta+\gamma$, we obtain
$$
s(t+1)-s(\alpha+1) \geq (\phi(y)-\alpha-1)s - {\phi(y)-\alpha-1 \choose 2}t = (\beta+\gamma)s - {\beta+\gamma \choose 2}t.
$$
Next note that the following inequalities
\[
(\beta+\gamma)s - {\beta+\gamma \choose 2}t \geq \beta s - {\beta \choose 2}t \quad \iff\quad \gamma \big(s- t\frac{2\beta+\gamma - 1}{2}\big) \geq 0,
\]
holds because $\gamma \geq 0$ and
\[
s\geq t(2\theta-1)\geq t(2\phi(y)-1) \geq t(2(\beta+\gamma+1)-1) \geq  t\frac{2\beta+\gamma - 1}{2}.
\]
Thus we may conclude
$$
s(t+1)-s(\alpha+1) \geq \beta s - {\beta \choose 2}t.
$$
Rearranging gives the desired inequality.
 \qed

We remark that the method of applying inclusion-exclusion principle is the same as was used in \cite{shilla} to prove Theorem \ref{shill}, and in essence they are the same result; however the statement of Theorem \ref{shill} requires the graph in question to be distance-regular, while here we are essentially applying it to a distance-regular graph with a clique removed.

\begin{lemma} \label{case2}
Suppose $s >t(2\theta-1)$ and (B) holds, then $s\leq (t+1)^2\theta\binom{\beta}{2}^{-1}$.
\end{lemma}

{\bf Proof:} We prove the result by establishing a lower bound on $|\tilde V|$ and then use it in conjunction with Lemma \ref{steps}. Consider $\DD(y)$ for any $y\in \tilde V$. If $x\in V(\DD(y))$ and $x\in V(G)\setminus \tilde V$ (thus $\phi(x)=t+1$), then by Lemma \ref{oneguy}, $x$ must belong to a clique of order $s$ in $\DD(y)$. By condition (B), we see that there must be at most $\alpha s$ such vertices in $\DD(y)$. Therefore,
\[
|\Ga(y) \cap \tilde V| \geq k-\alpha s = \beta s.
\]

Next we want to count the vertices $z\in \tilde V$ with $d(y,z)=2$. By applying the same argument as before, each $x\in \Ga(y) \cap \tilde V$ must also have at least $\beta s$ neighbours in $\tilde V$. Out of these, we must remove $y$ and those already in $\Ga(y)$. There are at most $|\{y\}\cup \Ga(y,x)| = s$ of those. Finally, for each $z\in \tilde V$ with $d(y,z)=2$, it is counted at most $|\Ga(y,z)|=t+1$ times via elements of $\Ga(y)$. Thus a (very weak) lower bound on $|\tilde V|$ is given by
\[
|\tilde V|\geq 1+\beta s+\frac{(\beta s) (\beta s-s)}{t+1} \geq \frac{(\beta s) (\beta s-s)}{t+1} = \binom{\beta}{2} \frac{2s^2}{t+1}.
\]

Combining this with Lemma \ref{steps}, we have
\[
\binom{\beta}{2} \frac{2s^2}{t+1} \leq 2(st+1)\theta.
\]
Using the fact that $st+1\leq s(t+1)$ and rearranging gives the required result. \qed


We are now in a position to finish the proof of Theorem \ref{newhotness}. The following summary will recap the key lemmas leading to the main result:

\begin{itemize}
\item Lemma \ref{lowerboundonphi} showed that $\phi(x)\geq t+1$ for all $x\in V(G)$.
\item Lemma \ref{itsgeo} showed that if $\phi(x)= t+1$ for all $x\in V(G)$, then we in fact have $GQ(s,t)$ and $s\leq t^2$.
\item Lemma \ref{bound1} showed that if $\phi(x)\geq \theta+1$ for some $x\in V(G)$, then $s \leq \frac{t}{\theta-t}\binom{\theta+1}{2}$.
\item So we are left with the case where $t+1\leq \phi(x)\leq \theta$ for all $x\in V(G)$ and $\phi(x)\geq t+2$ for some $x\in V(G)$. Let us henceforth assume that $s >t(2\theta-1)$. There are two remaining cases (A) and (B).
\item Lemma \ref{case1} showed that if (A) holds, then $s\leq \binom{\beta}{2} t$.
\item Lemma \ref{case2} showed that if (B) holds, then $s\leq (t+1)^2\theta\binom{\beta}{2}^{-1}$.
\end{itemize}
The combination of lemmas directly implies Theorem \ref{newhotness}, i.e., for any $PGQ(s,t)$,
\begin{gather}\label{eqbounds}
s \leq \max\left\{\frac{t}{\theta-t}\binom{\theta+1}{2},\ t(2\theta-1),\ \binom{\beta}{2}t,\ (t+1)^2\theta\binom{\beta}{2}^{-1}\right\}.
\end{gather}
\qed

{\bf Proof of Corollary \ref{iu}}
By considering eigenvalue multiplicities, a necessary condition for a $PQG(s,t)$ to exist is that $(s+t)$ must divide $s(s+1)t(t+1)$ (see \cite{godroy}). In the case of $t=2$ this means that $(s+2)$ divides $6s(s+1)$, which implies that $(s+2)$ divides $12$. Hence $s\leq 10$, which is less than the required bound of 12.

\vski

For $t\geq 3$,
choose $\theta=\left\lfloor \frac{4t}{3}+1 \right\rfloor$ and $\beta=\left\lceil 2\sqrt{t} \right\rceil$. We begin by showing that the third and fourth terms of \eqref{eqbounds} are bounded by $\frac{8}{3}t^2$. For the third term, since $2\sqrt{t}+1 \leq \frac{8}{3}\sqrt{t}$ for $t\geq 3$,
\[
\binom{\beta}{2}t \leq \frac{(2\sqrt{t}+1)(2\sqrt{t})}{2}t \leq \frac{8}{3} t^2.
\]
For the fourth term, since $\theta\leq \frac{4}{3}(t+1)$, $2\sqrt{t}-1\geq \frac{4}{3}\sqrt{t}$ and $(t+1)^3\leq \frac{8}{3}t^3$ for $t\geq 3$,
\[
(t+1)^2\theta\binom{\beta}{2}^{-1} \leq \frac{4}{3}(t+1)^3 \left( \frac{(2\sqrt{t})(2\sqrt{t}-1)}{2}\right)^{-1} \leq \frac{(t+1)^3 }{ t} \leq \frac{8}{3} t^2.
\]
For the second term of \eqref{eqbounds}, the required bound is implied by the following inequality
\begin{align}\label{eqboundsterm2}
t(2\theta-1)=t\left(2\left\lfloor \frac{4t}{3} \right\rfloor +1\right) \leq t\left\lfloor \frac{8t}{3}+1 \right\rfloor.
\end{align}

It remains to bound the first term of \eqref{eqbounds}. We will break it into three cases according to the residue of $t$ in modulo 3.
If $3|(t+1)$, then $\theta=\frac{4t+1}{3}$ and we have
\begin{align}\label{eqboundsterm1}
\frac{t}{\theta-t}\binom{\theta+1}{2} =\frac{t}{\frac{t+1}{3}}\frac{(\frac{4t+4}{3})(\frac{4t+1}{3})}{2}=\frac{t(8t+2)}{3}.
\end{align}
If $3|(t+2)$, then $\theta=\frac{4t+2}{3}$ and we have
\begin{align*}
\frac{t}{\theta-t}\binom{\theta+1}{2} =\frac{t}{\frac{t+2}{3}}\frac{(\frac{4t+5}{3})(\frac{4t+2}{3})}{2}\leq  \frac{t(8t+1)}{3} &\iff (4t+5)(4t+2)\leq 2(8t+1)(t+2) \\
&\iff 16t^2+28t+10 \leq 16t^2+34t+4.
\end{align*}
Finally, if $3|t$, then $\theta=\frac{4t+3}{3}$ and we have
\begin{align*}
\frac{t}{\theta-t}\binom{\theta+1}{2} =\frac{t}{\frac{t+3}{3}}\frac{(\frac{4t+6}{3})(\frac{4t+3}{3})}{2}\leq  \frac{t(8t+3)}{3} &\iff (4t+6)(4t+3)\leq 2(8t+3)(t+3) \\
&\iff 16t^2+36t+18 \leq 16t^2+54t+18.
\end{align*}
Therefore in all cases, the required inequality has been established and the proof is complete.
\qed

We remark that since \eqref{eqboundsterm2} achieves equality when $3|t$ or $3|(t+2)$, while \eqref{eqboundsterm1} holds when $3|(t+1)$, Corollary \ref{iu} in fact provides the tightest possible bound that can be derived from Theorem \ref{newhotness} for every $t$.

%
%

\section{Concluding remarks}

A computer search has uncovered the following intersection arrays with $t \leq 10$ which are ruled out by our Theorem \ref{newhotness} but which meet all other known feasibility conditions (Theorem \ref{oldbusted} as well as integral eigenvalue multiplicities).

$$
\begin{array}{rrrrrr}
  s & t & v & k & \la & \mu \\
  \hline

 56 & 4 & 12825 & 280 & 55 & 5 \\
 95 & 5 & 45696 & 570 & 94 & 6 \\
 120 & 6 & 87241 & 840 & 119 & 7 \\
 134 & 6 & 108675 & 938 & 133 & 7 \\
 140 & 7 & 138321 & 1120 & 139 & 8 \\
 161 & 7 & 182736 & 1288 & 160 & 8 \\
 189 & 7 & 251560 & 1512 & 188 & 8 \\
 184 & 8 & 272505 & 1656 & 183 & 9 \\
 216 & 8 & 375193 & 1944 & 215 & 9 \\
 244 & 8 & 478485 & 2196 & 243 & 9 \\
 280 & 8 & 629721 & 2520 & 279 & 9 \\
 328 & 8 & 863625 & 2952 & 327 & 9 \\
 231 & 9 & 482560 & 2310 & 230 & 10 \\
 261 & 9 & 615700 & 2610 & 260 & 10 \\
 315 & 9 & 896176 & 3150 & 314 & 10 \\
 351 & 9 & 1112320 & 3510 & 350 & 10 \\
 396 & 9 & 1415305 & 3960 & 395 & 10 \\
 423 & 9 & 1614592 & 4230 & 422 & 10 \\
 290 & 10 & 844191 & 3190 & 289 & 11 \\
 320 & 10 & 1027521 & 3520 & 319 & 11 \\
 386 & 10 & 1494207 & 4246 & 385 & 11 \\
 440 & 10 & 1940841 & 4840 & 439 & 11 \\
 485 & 10 & 2357586 & 5335 & 484 & 11 \\
 540 & 10 & 2921941 & 5940 & 539 & 11 \\
 650 & 10 & 4232151 & 7150 & 649 & 11 \\

\end{array}
$$

\section{Acknowledgements}

This work started at Monash University when J.H.K was visiting G.M. We thank Monash University for its hospitality. Also, it was (partially) done while J.P. was working at Wonkwang University, and it was worked on while J.H.K and G.M. were visiting Wonkwang University. We thank Wonkwang University for its hospitality. J.P. is supported by Basic Research Program through the National Research Foundation of Korea funded by Ministry of Education (NRF-2017R1D1A1B03032016). J.H.K. has been partially supported by the National Natural Science Foundation of China (Grants No. 11471009 and No. 11671376) and by the Anhui Initiative in Quantum Information Technologies (Grant No. AHY150200). G.M. has been partially supported by the Australian Research Council (Grants DP0988483 and DE140101201). I.G. has been partially supported by the Australian Research Council (Grant DP170101227).

\bibliographystyle{alpha}
\bibliography{Nonexistence}

\end{document}